\documentclass[runningheads,a4paper]{llncs}
\usepackage{amssymb}
\usepackage{graphicx}
\usepackage{url}
\newcommand{\keywords}[1]{\par\addvspace\baselineskip
\noindent\keywordname\enspace\ignorespaces#1}

\usepackage{amsmath}

\usepackage{verbatim}

\usepackage{tikz}
\definecolor{mediumspringgreen}{rgb}{0.0, 0.98039215, 0.60392156}

\newcommand{\bb}{\mathbb}
\newcommand{\R}{\bb R}

\newcommand{\Z}{\bb Z}

\newcommand{\B}{B}
\renewcommand{\P}{\mathcal{P}}

\newcommand\st{:}
\newcommand{\setcond}[2]{\left\{ #1 \,:\, #2 \right\}}

\DeclareMathOperator    \relint         {rel\,int}
\DeclareMathOperator    \verts          {vert}

\chardef\Myunderscore=`\_
\usepackage{hyperref}  
\newcommand\underscore{\Myunderscore\allowbreak}
\let\_=\underscore

\newcommand\githubsearchurl{https://github.com/mkoeppe/infinite-group-relaxation-code/search}
\usepackage{pgf}
\pgfkeyssetvalue{/sagefunc/bccz_counterexample}{\href{\githubsearchurl?q=\%22def+bccz_counterexample(\%22}{\sage{bccz\underscore{}counterexample}}}%)
\pgfkeyssetvalue{/sagefunc/bcdsp_arbitrary_slope}{\href{\githubsearchurl?q=\%22def+bcdsp_arbitrary_slope(\%22}{\sage{bcdsp\underscore{}arbitrary\underscore{}slope}}}%)
\pgfkeyssetvalue{/sagefunc/bhk_gmi_irrational}{\href{\githubsearchurl?q=\%22def+bhk_gmi_irrational(\%22}{\sage{bhk\underscore{}gmi\underscore{}irrational}}}%)
\pgfkeyssetvalue{/sagefunc/bhk_irrational}{\href{\githubsearchurl?q=\%22def+bhk_irrational(\%22}{\sage{bhk\underscore{}irrational}}}%)
\pgfkeyssetvalue{/sagefunc/bhk_slant_irrational}{\href{\githubsearchurl?q=\%22def+bhk_slant_irrational(\%22}{\sage{bhk\underscore{}slant\underscore{}irrational}}}%)
\pgfkeyssetvalue{/sagefunc/chen_4_slope}{\href{\githubsearchurl?q=\%22def+chen_4_slope(\%22}{\sage{chen\underscore{}4\underscore{}slope}}}%)
\pgfkeyssetvalue{/sagefunc/dg_2_step_mir}{\href{\githubsearchurl?q=\%22def+dg_2_step_mir(\%22}{\sage{dg\underscore{}2\underscore{}step\underscore{}mir}}}%)
\pgfkeyssetvalue{/sagefunc/dg_2_step_mir_limit}{\href{\githubsearchurl?q=\%22def+dg_2_step_mir_limit(\%22}{\sage{dg\underscore{}2\underscore{}step\underscore{}mir\underscore{}limit}}}%)
\pgfkeyssetvalue{/sagefunc/drlm_2_slope_limit}{\href{\githubsearchurl?q=\%22def+drlm_2_slope_limit(\%22}{\sage{drlm\underscore{}2\underscore{}slope\underscore{}limit}}}%)
\pgfkeyssetvalue{/sagefunc/drlm_3_slope_limit}{\href{\githubsearchurl?q=\%22def+drlm_3_slope_limit(\%22}{\sage{drlm\underscore{}3\underscore{}slope\underscore{}limit}}}%)
\pgfkeyssetvalue{/sagefunc/drlm_backward_3_slope}{\href{\githubsearchurl?q=\%22def+drlm_backward_3_slope(\%22}{\sage{drlm\underscore{}backward\underscore{}3\underscore{}slope}}}%)
\pgfkeyssetvalue{/sagefunc/gj_2_slope}{\href{\githubsearchurl?q=\%22def+gj_2_slope(\%22}{\sage{gj\underscore{}2\underscore{}slope}}}%)
\pgfkeyssetvalue{/sagefunc/gj_2_slope_repeat}{\href{\githubsearchurl?q=\%22def+gj_2_slope_repeat(\%22}{\sage{gj\underscore{}2\underscore{}slope\underscore{}repeat}}}%)
\pgfkeyssetvalue{/sagefunc/gj_forward_3_slope}{\href{\githubsearchurl?q=\%22def+gj_forward_3_slope(\%22}{\sage{gj\underscore{}forward\underscore{}3\underscore{}slope}}}%)
\pgfkeyssetvalue{/sagefunc/gmic}{\href{\githubsearchurl?q=\%22def+gmic(\%22}{\sage{gmic}}}%)
\pgfkeyssetvalue{/sagefunc/hildebrand_2_sided_discont_1_slope_1}{\href{\githubsearchurl?q=\%22def+hildebrand_2_sided_discont_1_slope_1(\%22}{\sage{hildebrand\underscore{}2\underscore{}sided\underscore{}discont\underscore{}1\underscore{}slope\underscore{}1}}}%)
\pgfkeyssetvalue{/sagefunc/hildebrand_2_sided_discont_2_slope_1}{\href{\githubsearchurl?q=\%22def+hildebrand_2_sided_discont_2_slope_1(\%22}{\sage{hildebrand\underscore{}2\underscore{}sided\underscore{}discont\underscore{}2\underscore{}slope\underscore{}1}}}%)
\pgfkeyssetvalue{/sagefunc/hildebrand_5_slope_22_1}{\href{\githubsearchurl?q=\%22def+hildebrand_5_slope_22_1(\%22}{\sage{hildebrand\underscore{}5\underscore{}slope\underscore{}22\underscore{}1}}}%)
\pgfkeyssetvalue{/sagefunc/hildebrand_5_slope_24_1}{\href{\githubsearchurl?q=\%22def+hildebrand_5_slope_24_1(\%22}{\sage{hildebrand\underscore{}5\underscore{}slope\underscore{}24\underscore{}1}}}%)
\pgfkeyssetvalue{/sagefunc/hildebrand_5_slope_28_1}{\href{\githubsearchurl?q=\%22def+hildebrand_5_slope_28_1(\%22}{\sage{hildebrand\underscore{}5\underscore{}slope\underscore{}28\underscore{}1}}}%)
\pgfkeyssetvalue{/sagefunc/hildebrand_discont_3_slope_1}{\href{\githubsearchurl?q=\%22def+hildebrand_discont_3_slope_1(\%22}{\sage{hildebrand\underscore{}discont\underscore{}3\underscore{}slope\underscore{}1}}}%)
\pgfkeyssetvalue{/sagefunc/kf_n_step_mir}{\href{\githubsearchurl?q=\%22def+kf_n_step_mir(\%22}{\sage{kf\underscore{}n\underscore{}step\underscore{}mir}}}%)
\pgfkeyssetvalue{/sagefunc/kzh_10_slope_1}{\href{\githubsearchurl?q=\%22def+kzh_10_slope_1(\%22}{\sage{kzh\underscore{}10\underscore{}slope\underscore{}1}}}%)
\pgfkeyssetvalue{/sagefunc/kzh_28_slope_1}{\href{\githubsearchurl?q=\%22def+kzh_28_slope_1(\%22}{\sage{kzh\underscore{}28\underscore{}slope\underscore{}1}}}%)
\pgfkeyssetvalue{/sagefunc/kzh_28_slope_2}{\href{\githubsearchurl?q=\%22def+kzh_28_slope_2(\%22}{\sage{kzh\underscore{}28\underscore{}slope\underscore{}2}}}%)
\pgfkeyssetvalue{/sagefunc/kzh_3_slope_param_extreme_1}{\href{\githubsearchurl?q=\%22def+kzh_3_slope_param_extreme_1(\%22}{\sage{kzh\underscore{}3\underscore{}slope\underscore{}param\underscore{}extreme\underscore{}1}}}%)
\pgfkeyssetvalue{/sagefunc/kzh_3_slope_param_extreme_2}{\href{\githubsearchurl?q=\%22def+kzh_3_slope_param_extreme_2(\%22}{\sage{kzh\underscore{}3\underscore{}slope\underscore{}param\underscore{}extreme\underscore{}2}}}%)
\pgfkeyssetvalue{/sagefunc/kzh_4_slope_param_extreme_1}{\href{\githubsearchurl?q=\%22def+kzh_4_slope_param_extreme_1(\%22}{\sage{kzh\underscore{}4\underscore{}slope\underscore{}param\underscore{}extreme\underscore{}1}}}%)
\pgfkeyssetvalue{/sagefunc/kzh_5_slope_fulldim_1}{\href{\githubsearchurl?q=\%22def+kzh_5_slope_fulldim_1(\%22}{\sage{kzh\underscore{}5\underscore{}slope\underscore{}fulldim\underscore{}1}}}%)
\pgfkeyssetvalue{/sagefunc/kzh_5_slope_fulldim_2}{\href{\githubsearchurl?q=\%22def+kzh_5_slope_fulldim_2(\%22}{\sage{kzh\underscore{}5\underscore{}slope\underscore{}fulldim\underscore{}2}}}%)
\pgfkeyssetvalue{/sagefunc/kzh_5_slope_fulldim_3}{\href{\githubsearchurl?q=\%22def+kzh_5_slope_fulldim_3(\%22}{\sage{kzh\underscore{}5\underscore{}slope\underscore{}fulldim\underscore{}3}}}%)
\pgfkeyssetvalue{/sagefunc/kzh_5_slope_fulldim_4}{\href{\githubsearchurl?q=\%22def+kzh_5_slope_fulldim_4(\%22}{\sage{kzh\underscore{}5\underscore{}slope\underscore{}fulldim\underscore{}4}}}%)
\pgfkeyssetvalue{/sagefunc/kzh_5_slope_fulldim_5}{\href{\githubsearchurl?q=\%22def+kzh_5_slope_fulldim_5(\%22}{\sage{kzh\underscore{}5\underscore{}slope\underscore{}fulldim\underscore{}5}}}%)
\pgfkeyssetvalue{/sagefunc/kzh_5_slope_fulldim_covers_1}{\href{\githubsearchurl?q=\%22def+kzh_5_slope_fulldim_covers_1(\%22}{\sage{kzh\underscore{}5\underscore{}slope\underscore{}fulldim\underscore{}covers\underscore{}1}}}%)
\pgfkeyssetvalue{/sagefunc/kzh_5_slope_fulldim_covers_2}{\href{\githubsearchurl?q=\%22def+kzh_5_slope_fulldim_covers_2(\%22}{\sage{kzh\underscore{}5\underscore{}slope\underscore{}fulldim\underscore{}covers\underscore{}2}}}%)
\pgfkeyssetvalue{/sagefunc/kzh_5_slope_fulldim_covers_3}{\href{\githubsearchurl?q=\%22def+kzh_5_slope_fulldim_covers_3(\%22}{\sage{kzh\underscore{}5\underscore{}slope\underscore{}fulldim\underscore{}covers\underscore{}3}}}%)
\pgfkeyssetvalue{/sagefunc/kzh_5_slope_fulldim_covers_4}{\href{\githubsearchurl?q=\%22def+kzh_5_slope_fulldim_covers_4(\%22}{\sage{kzh\underscore{}5\underscore{}slope\underscore{}fulldim\underscore{}covers\underscore{}4}}}%)
\pgfkeyssetvalue{/sagefunc/kzh_5_slope_fulldim_covers_5}{\href{\githubsearchurl?q=\%22def+kzh_5_slope_fulldim_covers_5(\%22}{\sage{kzh\underscore{}5\underscore{}slope\underscore{}fulldim\underscore{}covers\underscore{}5}}}%)
\pgfkeyssetvalue{/sagefunc/kzh_5_slope_fulldim_covers_6}{\href{\githubsearchurl?q=\%22def+kzh_5_slope_fulldim_covers_6(\%22}{\sage{kzh\underscore{}5\underscore{}slope\underscore{}fulldim\underscore{}covers\underscore{}6}}}%)
\pgfkeyssetvalue{/sagefunc/kzh_5_slope_q22_f10_1}{\href{\githubsearchurl?q=\%22def+kzh_5_slope_q22_f10_1(\%22}{\sage{kzh\underscore{}5\underscore{}slope\underscore{}q22\underscore{}f10\underscore{}1}}}%)
\pgfkeyssetvalue{/sagefunc/kzh_5_slope_q22_f10_2}{\href{\githubsearchurl?q=\%22def+kzh_5_slope_q22_f10_2(\%22}{\sage{kzh\underscore{}5\underscore{}slope\underscore{}q22\underscore{}f10\underscore{}2}}}%)
\pgfkeyssetvalue{/sagefunc/kzh_5_slope_q22_f10_3}{\href{\githubsearchurl?q=\%22def+kzh_5_slope_q22_f10_3(\%22}{\sage{kzh\underscore{}5\underscore{}slope\underscore{}q22\underscore{}f10\underscore{}3}}}%)
\pgfkeyssetvalue{/sagefunc/kzh_5_slope_q22_f10_4}{\href{\githubsearchurl?q=\%22def+kzh_5_slope_q22_f10_4(\%22}{\sage{kzh\underscore{}5\underscore{}slope\underscore{}q22\underscore{}f10\underscore{}4}}}%)
\pgfkeyssetvalue{/sagefunc/kzh_5_slope_q22_f2_1}{\href{\githubsearchurl?q=\%22def+kzh_5_slope_q22_f2_1(\%22}{\sage{kzh\underscore{}5\underscore{}slope\underscore{}q22\underscore{}f2\underscore{}1}}}%)
\pgfkeyssetvalue{/sagefunc/kzh_6_slope_1}{\href{\githubsearchurl?q=\%22def+kzh_6_slope_1(\%22}{\sage{kzh\underscore{}6\underscore{}slope\underscore{}1}}}%)
\pgfkeyssetvalue{/sagefunc/kzh_6_slope_fulldim_covers_1}{\href{\githubsearchurl?q=\%22def+kzh_6_slope_fulldim_covers_1(\%22}{\sage{kzh\underscore{}6\underscore{}slope\underscore{}fulldim\underscore{}covers\underscore{}1}}}%)
\pgfkeyssetvalue{/sagefunc/kzh_6_slope_fulldim_covers_2}{\href{\githubsearchurl?q=\%22def+kzh_6_slope_fulldim_covers_2(\%22}{\sage{kzh\underscore{}6\underscore{}slope\underscore{}fulldim\underscore{}covers\underscore{}2}}}%)
\pgfkeyssetvalue{/sagefunc/kzh_6_slope_fulldim_covers_3}{\href{\githubsearchurl?q=\%22def+kzh_6_slope_fulldim_covers_3(\%22}{\sage{kzh\underscore{}6\underscore{}slope\underscore{}fulldim\underscore{}covers\underscore{}3}}}%)
\pgfkeyssetvalue{/sagefunc/kzh_6_slope_fulldim_covers_4}{\href{\githubsearchurl?q=\%22def+kzh_6_slope_fulldim_covers_4(\%22}{\sage{kzh\underscore{}6\underscore{}slope\underscore{}fulldim\underscore{}covers\underscore{}4}}}%)
\pgfkeyssetvalue{/sagefunc/kzh_6_slope_fulldim_covers_5}{\href{\githubsearchurl?q=\%22def+kzh_6_slope_fulldim_covers_5(\%22}{\sage{kzh\underscore{}6\underscore{}slope\underscore{}fulldim\underscore{}covers\underscore{}5}}}%)
\pgfkeyssetvalue{/sagefunc/kzh_7_slope_1}{\href{\githubsearchurl?q=\%22def+kzh_7_slope_1(\%22}{\sage{kzh\underscore{}7\underscore{}slope\underscore{}1}}}%)
\pgfkeyssetvalue{/sagefunc/kzh_7_slope_2}{\href{\githubsearchurl?q=\%22def+kzh_7_slope_2(\%22}{\sage{kzh\underscore{}7\underscore{}slope\underscore{}2}}}%)
\pgfkeyssetvalue{/sagefunc/kzh_7_slope_3}{\href{\githubsearchurl?q=\%22def+kzh_7_slope_3(\%22}{\sage{kzh\underscore{}7\underscore{}slope\underscore{}3}}}%)
\pgfkeyssetvalue{/sagefunc/kzh_7_slope_4}{\href{\githubsearchurl?q=\%22def+kzh_7_slope_4(\%22}{\sage{kzh\underscore{}7\underscore{}slope\underscore{}4}}}%)
\pgfkeyssetvalue{/sagefunc/ll_strong_fractional}{\href{\githubsearchurl?q=\%22def+ll_strong_fractional(\%22}{\sage{ll\underscore{}strong\underscore{}fractional}}}%)
\pgfkeyssetvalue{/sagefunc/mlr_cpl3_a_2_slope}{\href{\githubsearchurl?q=\%22def+mlr_cpl3_a_2_slope(\%22}{\sage{mlr\underscore{}cpl3\underscore{}a\underscore{}2\underscore{}slope}}}%)
\pgfkeyssetvalue{/sagefunc/mlr_cpl3_b_3_slope}{\href{\githubsearchurl?q=\%22def+mlr_cpl3_b_3_slope(\%22}{\sage{mlr\underscore{}cpl3\underscore{}b\underscore{}3\underscore{}slope}}}%)
\pgfkeyssetvalue{/sagefunc/mlr_cpl3_c_3_slope}{\href{\githubsearchurl?q=\%22def+mlr_cpl3_c_3_slope(\%22}{\sage{mlr\underscore{}cpl3\underscore{}c\underscore{}3\underscore{}slope}}}%)
\pgfkeyssetvalue{/sagefunc/mlr_cpl3_d_3_slope}{\href{\githubsearchurl?q=\%22def+mlr_cpl3_d_3_slope(\%22}{\sage{mlr\underscore{}cpl3\underscore{}d\underscore{}3\underscore{}slope}}}%)
\pgfkeyssetvalue{/sagefunc/mlr_cpl3_f_2_or_3_slope}{\href{\githubsearchurl?q=\%22def+mlr_cpl3_f_2_or_3_slope(\%22}{\sage{mlr\underscore{}cpl3\underscore{}f\underscore{}2\underscore{}or\underscore{}3\underscore{}slope}}}%)
\pgfkeyssetvalue{/sagefunc/mlr_cpl3_g_3_slope}{\href{\githubsearchurl?q=\%22def+mlr_cpl3_g_3_slope(\%22}{\sage{mlr\underscore{}cpl3\underscore{}g\underscore{}3\underscore{}slope}}}%)
\pgfkeyssetvalue{/sagefunc/mlr_cpl3_h_2_slope}{\href{\githubsearchurl?q=\%22def+mlr_cpl3_h_2_slope(\%22}{\sage{mlr\underscore{}cpl3\underscore{}h\underscore{}2\underscore{}slope}}}%)
\pgfkeyssetvalue{/sagefunc/mlr_cpl3_k_2_slope}{\href{\githubsearchurl?q=\%22def+mlr_cpl3_k_2_slope(\%22}{\sage{mlr\underscore{}cpl3\underscore{}k\underscore{}2\underscore{}slope}}}%)
\pgfkeyssetvalue{/sagefunc/mlr_cpl3_l_2_slope}{\href{\githubsearchurl?q=\%22def+mlr_cpl3_l_2_slope(\%22}{\sage{mlr\underscore{}cpl3\underscore{}l\underscore{}2\underscore{}slope}}}%)
\pgfkeyssetvalue{/sagefunc/mlr_cpl3_n_3_slope}{\href{\githubsearchurl?q=\%22def+mlr_cpl3_n_3_slope(\%22}{\sage{mlr\underscore{}cpl3\underscore{}n\underscore{}3\underscore{}slope}}}%)
\pgfkeyssetvalue{/sagefunc/mlr_cpl3_o_2_slope}{\href{\githubsearchurl?q=\%22def+mlr_cpl3_o_2_slope(\%22}{\sage{mlr\underscore{}cpl3\underscore{}o\underscore{}2\underscore{}slope}}}%)
\pgfkeyssetvalue{/sagefunc/mlr_cpl3_p_2_slope}{\href{\githubsearchurl?q=\%22def+mlr_cpl3_p_2_slope(\%22}{\sage{mlr\underscore{}cpl3\underscore{}p\underscore{}2\underscore{}slope}}}%)
\pgfkeyssetvalue{/sagefunc/mlr_cpl3_q_2_slope}{\href{\githubsearchurl?q=\%22def+mlr_cpl3_q_2_slope(\%22}{\sage{mlr\underscore{}cpl3\underscore{}q\underscore{}2\underscore{}slope}}}%)
\pgfkeyssetvalue{/sagefunc/mlr_cpl3_r_2_slope}{\href{\githubsearchurl?q=\%22def+mlr_cpl3_r_2_slope(\%22}{\sage{mlr\underscore{}cpl3\underscore{}r\underscore{}2\underscore{}slope}}}%)
\pgfkeyssetvalue{/sagefunc/rlm_dpl1_extreme_3a}{\href{\githubsearchurl?q=\%22def+rlm_dpl1_extreme_3a(\%22}{\sage{rlm\underscore{}dpl1\underscore{}extreme\underscore{}3a}}}%)
\pgfkeyssetvalue{/sagefunc/automorphism}{\href{\githubsearchurl?q=\%22def+automorphism(\%22}{\sage{automorphism}}}%)
\pgfkeyssetvalue{/sagefunc/multiplicative_homomorphism}{\href{\githubsearchurl?q=\%22def+multiplicative_homomorphism(\%22}{\sage{multiplicative\underscore{}homomorphism}}}%)
\pgfkeyssetvalue{/sagefunc/projected_sequential_merge}{\href{\githubsearchurl?q=\%22def+projected_sequential_merge(\%22}{\sage{projected\underscore{}sequential\underscore{}merge}}}%)
\pgfkeyssetvalue{/sagefunc/restrict_to_finite_group}{\href{\githubsearchurl?q=\%22def+restrict_to_finite_group(\%22}{\sage{restrict\underscore{}to\underscore{}finite\underscore{}group}}}%)
\pgfkeyssetvalue{/sagefunc/interpolate_to_infinite_group}{\href{\githubsearchurl?q=\%22def+interpolate_to_infinite_group(\%22}{\sage{interpolate\underscore{}to\underscore{}infinite\underscore{}group}}}%)
\pgfkeyssetvalue{/sagefunc/two_slope_fill_in}{\href{\githubsearchurl?q=\%22def+two_slope_fill_in(\%22}{\sage{two\underscore{}slope\underscore{}fill\underscore{}in}}}%)
\pgfkeyssetvalue{/sagefunc/generate_example_e_for_psi_n}{\href{\githubsearchurl?q=\%22def+generate_example_e_for_psi_n(\%22}{\sage{generate\underscore{}example\underscore{}e\underscore{}for\underscore{}psi\underscore{}n}}}%)
\pgfkeyssetvalue{/sagefunc/chen_3_slope_not_extreme}{\href{\githubsearchurl?q=\%22def+chen_3_slope_not_extreme(\%22}{\sage{chen\underscore{}3\underscore{}slope\underscore{}not\underscore{}extreme}}}%)
\pgfkeyssetvalue{/sagefunc/psi_n_in_bccz_counterexample_construction}{\href{\githubsearchurl?q=\%22def+psi_n_in_bccz_counterexample_construction(\%22}{\sage{psi\underscore{}n\underscore{}in\underscore{}bccz\underscore{}counterexample\underscore{}construction}}}%)
\pgfkeyssetvalue{/sagefunc/gomory_fractional}{\href{\githubsearchurl?q=\%22def+gomory_fractional(\%22}{\sage{gomory\underscore{}fractional}}}%)
\pgfkeyssetvalue{/sagefunc/not_minimal_2}{\href{\githubsearchurl?q=\%22def+not_minimal_2(\%22}{\sage{not\underscore{}minimal\underscore{}2}}}%)
\pgfkeyssetvalue{/sagefunc/not_extreme_1}{\href{\githubsearchurl?q=\%22def+not_extreme_1(\%22}{\sage{not\underscore{}extreme\underscore{}1}}}%)
\pgfkeyssetvalue{/sagefunc/kzh_2q_example_1}{\href{\githubsearchurl?q=\%22def+kzh_2q_example_1(\%22}{\sage{kzh\underscore{}2q\underscore{}example\underscore{}1}}}%)
\pgfkeyssetvalue{/sagefunc/zhou_two_sided_discontinuous_cannot_assume_any_continuity}{\href{\githubsearchurl?q=\%22def+zhou_two_sided_discontinuous_cannot_assume_any_continuity(\%22}{\sage{zhou\underscore{}two\underscore{}sided\underscore{}discontinuous\underscore{}cannot\underscore{}assume\underscore{}any\underscore{}continuity}}}%)
\pgfkeyssetvalue{/sagefunc/extremality_test}{\href{\githubsearchurl?q=\%22def+extremality_test(\%22}{\sage{extremality\underscore{}test}}}%)
\pgfkeyssetvalue{/sagefunc/plot_2d_diagram}{\href{\githubsearchurl?q=\%22def+plot_2d_diagram(\%22}{\sage{plot\underscore{}2d\underscore{}diagram}}}%)
\pgfkeyssetvalue{/sagefunc/nice_field_values}{\href{\githubsearchurl?q=\%22def+nice_field_values(\%22}{\sage{nice\underscore{}field\underscore{}values}}}%)
\pgfkeyssetvalue{/sagefunc/ParametricRealFieldElement}{\href{\githubsearchurl?q=\%22def+ParametricRealFieldElement(\%22}{\sage{ParametricRealFieldElement}}}%)
\pgfkeyssetvalue{/sagefunc/ParametricRealField}{\href{\githubsearchurl?q=\%22def+ParametricRealField(\%22}{\sage{ParametricRealField}}}%)

\DeclareRobustCommand\sage[1]{\texttt{#1}}
\DeclareRobustCommand\sagefunc[1]{\pgfkeys{/sagefunc/#1}}

\makeatletter %% Fix for spacing of three stars
\def\@fnsymbol#1{\ensuremath{\ifcase#1\or\star\or{\star\star}\or
   {\star{\star}\star}\or \dagger\or \ddagger\or
   \mathchar "278\or \mathchar "27B\or \|\or **\or \dagger\dagger
   \or \ddagger\ddagger \else\@ctrerr\fi}}
\makeatother

\begin{document}

\mainmatter
%================================================

%==== FILL IN ====================================
\title{Software for cut-generating functions in the Gomory--Johnson
model and beyond%
\thanks{The authors gratefully acknowledge partial support from the National Science
  Foundation through grant DMS-1320051 awarded to M.~K\"oppe.}
}  % Full title
%% \subtitle{Power tools for the modern\\ cutgeneratingfunctionologist}
\titlerunning{Software for cut-generating functions} % Short title
\author{Chun Yu Hong\inst{1}\thanks{The first author's contribution was done
    during a Research Experiences for Undergraduates at the University of
    California, Davis. He was partially supported by the National Science
    Foundation through grant DMS-0636297 (VIGRE).
  } 
  \and Matthias K\"oppe\inst{2}\thanks{Corresponding author.} \and Yuan Zhou\inst{2}}
\authorrunning{Hong, K\"{o}ppe, Zhou}
\institute{
University of California, Berkeley, Department of Statistics, USA
\email{jcyhong@berkeley.edu},\\ 
\url{http://statistics.berkeley.edu/people/chun-yu-hong}
\and
University of California, Davis, Department of Mathematics, USA\\
\email{\texttt{\{mkoeppe,yzh\}@math.ucdavis.edu}},\\
\url{http://www.math.ucdavis.edu/~{mkoeppe,yzh}}\\ 
}
\maketitle

\begin{abstract}
  We present software for investigations with cut generating functions in the
  Gomory--Johnson model and extensions, implemented in the computer algebra system
  SageMath.
\end{abstract}

%------------------------------------------------------------
\section{Introduction}

Consider the following question from the theory of linear inequalities over the
reals: Given a (finite) system $Ax \leq b$, exactly which linear inequalities
$\langle a,x\rangle \leq \beta$ are \emph{valid}, i.e., satisfied for every
$x$ that satisfies the given system?  The answer is given, of course, by the
Farkas Lemma, or, equivalently, by the strong duality theory of linear
optimization.  As is well-known, this duality theory is symmetric: The dual of
a linear optimization problem is again a linear optimization problem, and the
dual of the dual is the original (primal) optimization problem.

The question becomes much harder when all or some of the variables are
constrained to be integers.  The theory of valid linear inequalities here is
called \emph{cutting plane theory}.  Over the past 60 years, a vast body of
research has been carried out on this topic, the largest part of it regarding
the polyhedral combinatorics of integer hulls of particular families of
problems.  The general theory again is equivalent to the duality theory of
integer linear optimization problems.  Here the dual objects are not linear,
but \emph{superadditive} (or subadditive) functionals, making the general form of
this theory infinite-dimensional even though the original problem started out
with only finitely many variables.

These superadditive (or subadditive) functionals appear in integer linear
optimization in various concrete forms, for example in the form of
\emph{dual-feasible functions}
\cite{alves-clautiaux-valerio-rietz-2016:dual-feasible-book}, 
\emph{superadditive lifting functions} 
\cite{louveaux-wolsey-2003:lifting-superadditivity-mir-survey}, 
and 
\emph{cut-generating functions}~\cite{conforti2013cut}. 

In the present paper, we describe some aspects of our software
\cite{infinite-group-relaxation-code} for
cut-generating functions in the classic 1-row Gomory--Johnson
\cite{infinite,infinite2} model.  In this theory, the main objects are the
so-called \emph{minimal valid functions}, which are the $\Z$-periodic,
subadditive functions $\pi\colon \R\to\R_+$ with $\pi(0)=0$, $\pi(f)=1$, that
satisfy the \emph{symmetry condition} $\pi(x) + \pi(f - x) = 1$ for all
$x\in\R$.  (Here $f$ is a fixed number.)
%% TODO: Refer to an example figure.
We refer the reader to the recent survey \cite{igp_survey,igp_survey_part_2}. 

Our software is a tool that enables mathematical exploration and research in
this domain. It can also be used in an educational setting, where it enables
hands-on teaching about modern cutting plane theory based on cut-generating
functions.  It removes the limitations of hand-written proofs, which would be
dominated by tedious case analysis.

%% Can be understood as the duality theory of the infinite-dimensional
%% finite-support integer optimization problems 

The first version of our software \cite{infinite-group-relaxation-code} was
written by the first author, C.~Y. Hong, during a Research Experience for
Undergraduates % at the University of California, Davis,
in summer 2013.  It was
later revised and extended by % the second author, 
M.~K\"oppe %, 
and again by % the
% third author,
Y.~Zhou.  The latter added an electronic compendium
\cite{zhou:extreme-notes} of extreme
functions found in the literature% , which has been described in the
% article~
, and added code that handles the case of
discontinuous functions.  Version 0.9 of our software was released in 2014 to
accompany the survey \cite{igp_survey,igp_survey_part_2}; the software has
received continuous updates by the second and third authors since.\footnote{Two
further undergraduate students contributed to our software.  P.~Xiao
contributed some documentation and tests.  M.~Sugiyama contributed additional
functions to the compendium, and added code for % to convert between group functions
% and
superadditive lifting functions.}

%%% \section{About SageMath as a research and education platform in optimization}

Our software is written in Python, making use of the convenient framework of
the open-source computer algebra system SageMath \cite{sage}.  It can be run
on a local installation of SageMath, or online via \emph{SageMathCloud}.

\section{Continuous and discontinuous piecewise linear $\Z$-periodic functions}

The main objects of our code are the $\Z$-periodic functions~$\pi\colon \R\to\R$.  Our code
is limited to the case of piecewise linear functions, which are allowed to be
discontinuous; see the definition below.  
In the following, we connect to the systematic notation introduced in
\cite[section 2.1]{basu-hildebrand-koeppe:equivariant}; see also % , which generalizes to the multi-row
% case (see the
% survey~
\cite{igp_survey,igp_survey_part_2}.
In our code, the periodicity of the functions is implicit; the functions
are represented by their restriction to the interval $[0,1]$.\footnote{The
  functions are instances of the class \sage{FastPiecewise}, which extends an existing
  SageMath class for piecewise linear functions.}
They can be constructed in various ways using Python functions named
\sage{piecewise\_function\_from\_breakpoints\_and\_values} etc.; see %   We suppress
% the details; 
the source code of the electronic compendium %provides many
for
examples.  We also suppress the details of the internal representation; instead we explain the
main ways in which the data of the function are accessed.

\begin{description}
\item[\sage{$\pi$.end\_points()}] is a list
  $0=x_0 < x_1 < \dots < x_{n-1} < x_n=1$ of possible
  breakpoints
  \begin{comment}
    \footnote{If the function~$\pi$ has been constructed with
      \sage{merge=True} (the default), then it is guaranteed that all end
      points $x_i$, with the possible exception of $0$ and $1$, are actual
      breakpoints of the $\Z$-periodic function~$\pi$.}
  \end{comment}
  of the function in $[0,1]$.  
  In the notation from \cite{basu-hildebrand-koeppe:equivariant,igp_survey,igp_survey_part_2}, 
  these endpoints are extended periodically as
  \begin{math}
    \B = \{\, x_0 + t, x_1 + t, \dots, x_{n-1}+t\st
    t\in\Z\,\}
  \end{math}.  Then the set of 0-dimensional faces is defined to be the collection of
  singletons, $\bigl\{\, \{ x \} \st x\in B\,\bigr\}$, 
  and the set of one-dimensional faces to be the collection of closed intervals,
  \begin{math}
    \bigl\{\, [x_i+t, x_{i+1}+t] \st i=0, \dots, {n-1} \text{ and } t\in\Z \,\bigr\}. 
  \end{math}
  Together, we obtain $\P = \P_{\B}$, a locally finite % one-dimensional
  polyhedral
  complex, % that is
  periodic modulo~$\Z$.
\item[\sage{$\pi$.values\_at\_end\_points()}] is a list 
  of the function values $\pi(x_i)$, %for
  $i=0, \dots,n$. 
  This list is most useful for continuous piecewise linear functions, as 
  indicated by \sage{$\pi$.is\_continuous()}, in which case the function is
  defined on the intervals $[x_i, x_{i+1}]$ by linear interpolation.
\item[\sage{$\pi$.limits\_at\_end\_points()}] provides data for the
  general, possibly discontinuous case in the form of 
  a list \sage{limits} of 3-tuples, 
  with
  \begin{align*}
    \sage{limits[$i$][0]} &= \pi(x_i) \\
    \sage{limits[$i$][1]} &= \pi(x_i^+) = \lim_{x\to x_i, x>x_i} \pi(x)\\
    \sage{limits[$i$][-1]} &= \pi(x_i^-) = \lim_{x\to x_i, x<x_i} \pi(x). 
  \end{align*}
  The function is defined on the open intervals $(x_i, x_{i+1})$ by linear
  interpolation of the limit values $\pi(x_i^+)$, $\pi(x_{i+1}^-)$. 
\item[\sage{$\pi$($x$)} and \sage{$\pi$.limits($x$)}] evaluate the function at
  $x$ and provide the 3-tuple of its limits at $x$, respectively.
\item[\sage{$\pi$.which\_function($x$)}] returns a linear function, denoted
  $\pi_I\colon\R\to\R$ in \cite{basu-hildebrand-koeppe:equivariant,igp_survey,igp_survey_part_2},
  where $I$ is the smallest face of~$\P$ containing $x$, so $\pi(x) =
  \pi_I(x)$ for $x \in \relint(I)$. 
\end{description}
%% 1) as a list of closed, open, or half-open intervals and a linear function 
Functions can be plotted using the standard SageMath function
\sage{plot($\pi$)}, or using our function
\sage{plot\_with\_colored\_slopes($\pi$)}, which assigns a different color to
each different slope value that a linear piece takes.\footnote{See also our
  function \sage{number\_of\_slopes}.  We refer the reader to \cite[section 2.4]{igp_survey} 
  for a discussion of the number of slopes of extreme functions, 
  and \cite{bcdsp:arbitrary-slopes} and \sagefunc{bcdsp_arbitrary_slope} for the
  latest developments in this direction.}
Examples of such functions are shown in Figures
\ref{fig:2d_diagram_with_cones} and \ref{fig:2d_diagrams_discontinuous_function}.

%%  and \sage{plot\_covered\_intervals}

\section{The diagrams of the decorated 2-dimensional polyhedral complex $\Delta\P$}

We now describe certain 2-dimensional diagrams which record the
subadditivity and additivity properties of a given function.  
These diagrams, in the continuous case, have appeared extensively in
\cite{igp_survey,igp_survey_part_2,zhou:extreme-notes}.  An example for the
discontinuous case appeared in \cite{zhou:extreme-notes}.  
We have engineered these diagrams from earlier forms that can be found in
\cite{tspace} (for the discussion of the \sage{merit\_index}) and in
\cite{basu-hildebrand-koeppe:equivariant}, to become 
power tools for the modern cutgeneratingfunctionologist.
Not only is the minimality of a given function immediately apparent on the
diagram, but also the extremality proof for a given class of piecewise minimal
valid functions follows a standard pattern that draws from these diagrams.  See
\cite[prelude]{igp_survey_part_2} and \cite[sections 2 and
4]{zhou:extreme-notes} for examples of such proofs.

\subsection{The polyhedral complex and its faces}
\begin{figure}[tp]
  \centering
  \includegraphics[width=.5\linewidth]{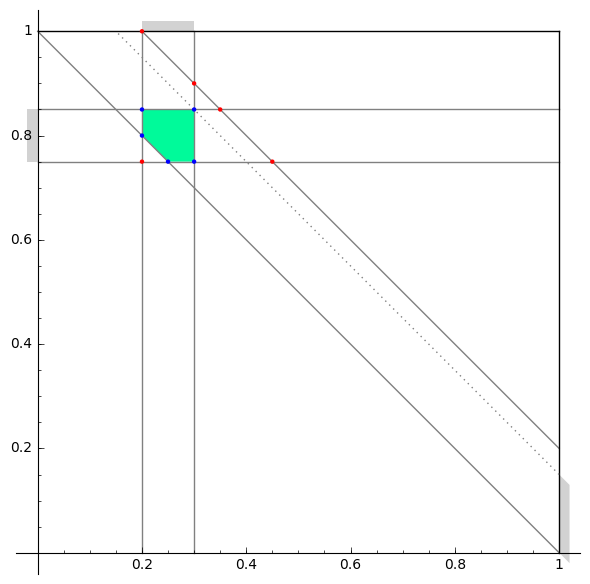}
  \caption{An example of a face $F = F(I, J, K)$ of the 2-dimensional
    polyhedral complex $\Delta\P$, set up by \sage{F = Face([[0.2, 0.3],
      [0.75, 0.85], [1, 1.2]])}.  It has vertices (\emph{blue})
    $(0.2, 0.85), (0.3, 0.75), (0.3, 0.85), (0.2, 0.8), (0.25, 0.75)$,
    whereas the other basic solutions (\emph{red})
    $(0.2, 0.75), (0.2, 1), (0.3, 0.9), (0.35, 0.85), (0.45, 0.75)$
    are filtered out because they are infeasible. 
    The face $F$ has projections (\emph{gray shadows})
    $I' = p_1(F) = [0.2, 0.3]$ (\emph{top border}), $J' = p_2(F) = [0.75,
    0.85]$ (\emph{left border}), and $K'
      = p_3(F) = [1, 1.15]$ (\emph{right border}). Note that $K'\subsetneq K$. 
  } 
  \label{fig:construct_a_face}
\end{figure}

Following
\cite{basu-hildebrand-koeppe:equivariant,igp_survey,igp_survey_part_2}, we
introduce the function 
\[\Delta\pi \colon \R \times \R \to \R,\quad \Delta\pi(x,y) =
  \pi(x)+\pi(y)-\pi(x+y),\] which measures the
slack in the subadditivity condition.\footnote{It is available in the code as
  \sage{delta\_pi($\pi$, $x$, $y$)}; in \cite{infinite}, it was called $\nabla(x,y)$.}  
Thus, if $\Delta\pi(x,y)<0$, subadditivity is violated at $(x, y)$; 
if $\Delta\pi(x,y)=0$, additivity holds at $(x,y)$; 
and if $\Delta\pi(x,y)>0$, we have strict subadditivity at $(x,y)$.
% It should be clear that
The piecewise linearity of $\pi(x)$ 
induces piecewise linearity of $\Delta\pi(x,y)$.  To express the domains of
linearity of $\Delta\pi(x,y)$, and thus domains of additivity and strict
subadditivity, we introduce the two-dimensional polyhedral complex
$\Delta\P$. 
The faces $F$ of the complex are defined as follows. Let $I, J, K \in
\mathcal{P}$, so each of $I, J, K$ is either a breakpoint of $\pi$ or a closed
interval delimited by two consecutive breakpoints. Then 
$$ F = F(I,J,K) = \setcond{\,(x,y) \in \R \times \R}{x \in I,\, y \in J,\, x + y \in
  K\,}.$$ 
In our code, a face is represented by an instance of the class \sage{Face}. 
It is constructed from $I, J, K$ and is represented by 
the list of vertices of $F$ and its projections $I'=p_1(F)$, $J'=p_2(F)$, $K'=
p_3(F)$,
where $p_1, p_2, p_3 \colon \R \times \R \to \R$ are defined as 
$p_1(x,y)=x$,  $p_2(x,y)=y$, $p_3(x,y) = x+y$.
%$I' := p_1(F), J'=p_2(F), K':=p_3(F)$ of $F$. (\tred{Define projections.}) 
The vertices $\verts(F)$ are obtained by first listing the basic solutions
$(x,y)$ where $x$, $y$, and $x+y$ are fixed to endpoints of $I$, $J$, and $K$,
respectively, and then filtering the feasible solutions. 
%% Simple formulas for the projections of $F$ are available ():
 %% \begin{alignat*}{2}
 %%    I' := p_1(F(I,J,K)) &= (K + (-J)) \cap I && \subseteq I, \\
 %%    J' := p_2(F(I,J,K)) &= (K+ (-I)) \cap J && \subseteq J, \\
 %%    K' := p_3(F(I,J,K)) &= (I+J) \cap K && \subseteq K.
 %%  \end{alignat*}
%%\tred{Note:} The code does not use these formulas to compute $I', J',
%%K'$. It computes 
The three projections are then computed from the list of
vertices.
\begin{comment}
  \footnote{We do not use the formulas for the projections given by
    \cite[Proposition 3.3]{bhk-IPCOext}, \cite[equation
    (3.11)]{igp_survey}.} %% of $\verts(F)$ and takes $\min$ and
  %% %% $\max$.
\end{comment}
Due to the $\Z$-periodicity of $\pi$, we can represent a face as a subset
of $[0,1]\times[0,1]$.  See \autoref{fig:construct_a_face} for an example. 
Because of the importance of the projection
$p_3(x,y)=x + y$, it is convenient to imagine a third, $(x+y)$-axis in
addition to the $x$-axis and the $y$-axis, which
traces the bottom border for $0 \leq x+y \leq 1$ and then the right border for
$1 \leq x+y \leq 2$.  To make room for this new axis, the $x$-axis should be
drawn on the top border of the diagram.

\subsection{\sage{plot\_2d\_diagram\_with\_cones}}

\begin{figure}[t]
\centering
\begin{minipage}{.49\textwidth}
\centering
\includegraphics[width=.8\linewidth]{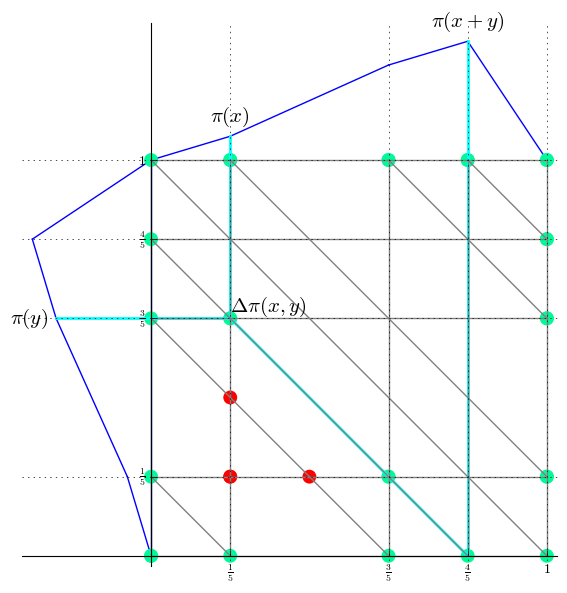}
\end{minipage}
\begin{minipage}{.49\textwidth}
\centering
\includegraphics[width=.8\linewidth]{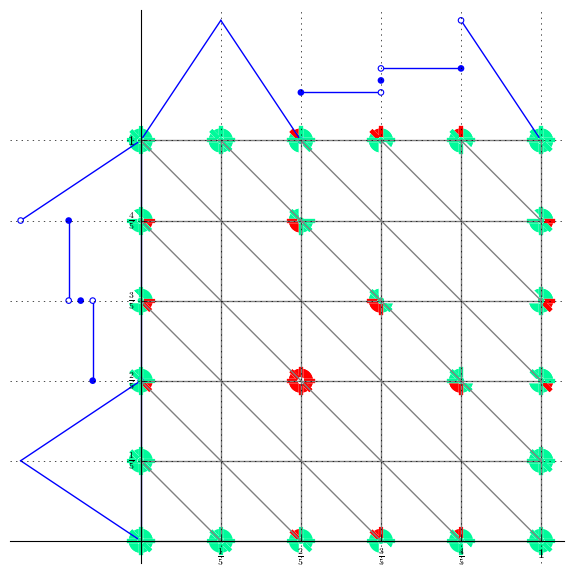}
\end{minipage}
\caption{Two diagrams of functions and their polyhedral complexes $\Delta\mathcal{P}$ with colored cones at $\verts(\Delta\mathcal{P})$, as plotted by the command \sage{plot\_2d\_diagram\_with\_cones(h)}. \textit{Left}, continuous function \sage{h = not\_minimal\_2()}. \textit{Right}, a random discontinuous function generated by \sage{h = random\_piecewise\_function(xgrid=5, ygrid=5, continuous\_proba=1/3, symmetry=True)}.}
\label{fig:2d_diagram_with_cones}
\end{figure}

We % are
now % ready to
explain the first version of the % complete 
2-dimensional
diagrams, plotted by the function
\sage{plot\_2d\_diagram\_with\_cones($\pi$)}; see
\autoref{fig:2d_diagram_with_cones}% for two examples
.
%% We are interested in the strict-subadditivity, additivity or non-subadditivity of the function $\pi$ over the faces of its polyhedral complex $\Delta\mathcal{P}$. All information can be visualized on the diagram plotted by the command \sage{plot\_2d\_diagram}. 
%% The diagram is in a two-dimensional space, with $3$ axes: 
%% \begin{itemize}
%% \item $x$-axis $(x \in [0,1])$ at the top border;
%% \item $y$-axis $(y \in [0,1])$ at the left border;
%% \item $x+y$-axis $(x+y \in [0,2])$ crooked at $x+y=1$ at the bottom and the right borders.
%% \end{itemize}
At the border of these diagrams, the function $\pi$ is shown twice
(\emph{blue}), 
along the $x$-axis (\emph{top border}) and along the $y$-axis (\emph{left
  border}). 
The solid grid lines in the diagrams are determined by the breakpoints of
$\pi$: vertical, horizontal and diagonal grid lines correspond to values where
$x$, $y$ and $x+y$ are breakpoints of $\pi$, respectively. 
The vertices of the complex $\Delta\mathcal{P}$ are the intersections of these
grid lines. %% Thus, there are three types of vertices in $\verts(\Delta\mathcal{P})$: 
%% \[(x,y), \quad (x, z-x) \quad \text{ and } \quad (z-y, y),\] where $x, y,
%% z$ are breakpoints of $\pi$, $0 \leq x \leq 1$, $0 \leq y \leq 1$ and $0
%% \leq z \leq 2$.  ----- OMITTED, repetitive

\textbf{In the continuous case}, % we investigate the function \[\Delta\pi
%   \colon \R \times \R \to \R, \Delta\pi(x,y) = \pi(x)+\pi(y)-\pi(x+y)\] at
% every vertices $(x, y)$, which measures the slack in the subadditivity
% constraint. W
we indicate the sign of $\Delta\pi(x,y)$ for all vertices by colored dots on
the diagram: \emph{red} indicates $\Delta\pi(x,y)<0$ (subadditivity is
violated); \emph{green} indicates $\Delta\pi(x,y)=0$ (additivity holds). 
%% \begin{itemize}
%% \item If $\Delta\pi(x,y)<0$, i.e., the subadditivity is violated at $(x, y)$, then the vertex $(x,y)$ is marked as a red dot on the diagram. 
%% \item If $\Delta\pi(x,y)=0$, i.e., the additivity holds at $(x,y)$, then the vertex $(x, y)$ is marked green.
%% \item If $\Delta\pi(x,y)>0$, i.e., the subadditivity is strict at $(x,y)$, then the vertex $(x,y)$ is marked white.
%% \end{itemize}
%%%% Note that if the first case happens at a vertex $(x, y)$, then the function
%%%% $\pi$ is non-subadditive and thus non-minimal.   
%%%% %% Since the diagram is symmetric with respect to $x=y$ and is periodic modulo $\Z \times \Z$, it suffices to consider the first two types of vertices $(x, y)$ and $(x, z-x)$, with $0 \leq x \leq y \leq 1$ and $x \leq z \leq x+1$; the other half of vertices can be obtained by swapping $x$ and $y$ values. 

\begin{example} 
  In \autoref{fig:2d_diagram_with_cones} (left), the vertex $(x, y)=(\frac{1}{5}, \frac{3}{5})$ is marked green, since 
\begin{align*}
\Delta\pi(\tfrac{1}{5},\tfrac{3}{5}) & =
                                       \pi(\tfrac{1}{5})+\pi(\tfrac{3}{5})-\pi(\tfrac{4}{5})
  %\\
%%& \quad \text{(read off the values of $\pi$ from the functions plotted at the borders of the diagram)} \\
%& 
= \tfrac{1}{5} +\tfrac{4}{5} -1 = 0.
\end{align*}
%% By the invariance of $\Delta\pi(x,y)$ under swapping $x$ and $y$,
%% $\Delta\pi(\frac{3}{5}, \frac{1}{5}) = \Delta\pi(\frac{1}{5},\frac{3}{5}) =0$,
%% hence the vertex $(\frac{3}{5}, \frac{1}{5})$ is also marked green.
\end{example} 

\textbf{In the discontinuous case}, beside the subadditivity slack $\Delta\pi(x,y)$ at a vertex $(x, y)$, one also needs to study the limit value of $\Delta\pi$ at the vertex $(x,y)$ approaching from the interior of a face $F \in \Delta\mathcal{P}$ containing the vertex $(x,y)$. This limit value is defined by
\[\Delta\pi_F(x,y) = \lim_{\substack{(u,v) \to (x,y)\\ (u,v) \in \relint(F)}}
  \Delta\pi(u,v), \quad \text{where } F \in \Delta\mathcal{P} \text{ such that
  } (x, y) \in F.\]  We indicate the sign of $\Delta\pi_F(x,y)$ by a colored
cone inside $F$ pointed at the vertex $(x, y)$ on the diagram. There could be
up to $12$ such cones (including rays for one-dimensional $F$) around a vertex $(x,
y)$.

\begin{example}
  In \autoref{fig:2d_diagram_with_cones} (right), the lower
  right corner $(x, y)=(\frac{2}{5}, \frac{4}{5})$ of the face
  $F = F(I, J, K)$ with
  $I =[\frac{1}{5}, \frac{2}{5}]$, $J=[\frac{4}{5}, 1]$, $K=[1, \frac{6}{5}]$ is
  green, since
  \begin{align*}
    \Delta\pi_F(x,y) & = \lim_{\substack{(u,v) \to (\frac{2}{5},\frac{4}{5})\\ (u,v) \in \relint(F)}} \Delta\pi(u,v) \\
                     & = \lim_{u\to\frac{2}{5},\; u < \frac{2}{5}}\pi(u)+\lim_{v\to\frac{4}{5}, \; v > \frac{4}{5}}\pi(v)-\lim_{w\to\frac{6}{5},\; w <\frac{6}{5}}\pi(w) \\
                     & = \pi(\tfrac{2}{5}^-) +\pi(\tfrac{4}{5}^+) -\pi(\tfrac{1}{5}^-) \quad \text{(as } \pi(\tfrac{6}{5}^-) =\pi(\tfrac{1}{5}^-) \text{ by periodicity)} \\
                     & = 0 + 1 - 1 = 0.
  \end{align*}
  The horizontal ray to the left of the same vertex
  $(x, y)=(\frac{2}{5}, \frac{4}{5})$ is red, because approaching from the
  one-dimensional face $F' = F(I', J', K')$ that contains $(x,y)$, with
  $I' =[\frac{1}{5}, \frac{2}{5}]$, $J'=\{\frac{4}{5}\}$, $K'=[1, \frac{6}{5}]$,
  we have the limit value
  \[\Delta\pi_{F'}(x,y) = \lim_{\substack{(u,v) \to
        (\frac{2}{5},\frac{4}{5})\\ (u,v) \in \relint(F')}} \hspace{-1.5em}\Delta\pi(u,v) =
    \lim_{\substack{u\to\frac{2}{5} \\ u <
        \frac{2}{5}}}\pi(u)+\pi(\tfrac{4}{5})-\lim_{\substack{w\to\frac{6}{5}\\
        w <\frac{6}{5}}}\pi(w) = 0 + \tfrac{3}{5} - 1 < 0.\]
  % The command \sage{plot\_2d\_diagram\_with\_cones()} draws the color-coded
  % cones.
\end{example}

\subsection{\sage{plot\_2d\_diagram} and additive faces}

\begin{figure}[tp]
\centering
\begin{minipage}{.49\textwidth}
\centering
\includegraphics[width=.8\linewidth]{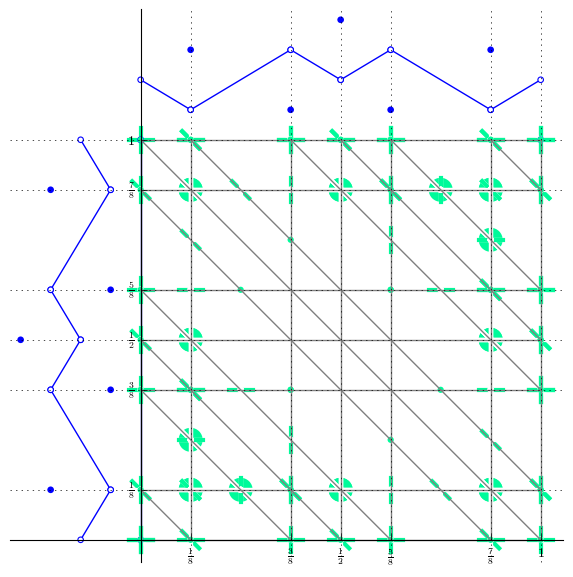}
\end{minipage}
\begin{minipage}{.49\textwidth}
\centering
\includegraphics[width=.8\linewidth]{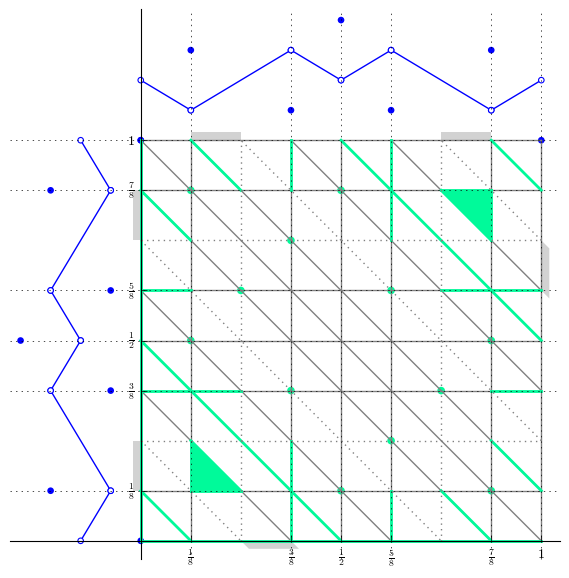}
\end{minipage}
\caption{Diagrams of $\Delta\mathcal{P}$ of a discontinuous function \sage{h = hildebrand\underscore 2\underscore sided\underscore discont\underscore 2\underscore slope\underscore 1()}, with (\textit{left}) additive limiting cones as plotted by the command \sage{plot\_2d\_diagram\_with\_cones(h)}; (\textit{right}) additive faces as plotted by the command \sage{plot\_2d\_diagram(h)}.}
\label{fig:2d_diagrams_discontinuous_function}
\end{figure}

Now assume that $\pi$ is a subadditive function. Then there are no red dots or cones on the above diagram of the complex $\Delta\mathcal{P}$. %% Next, we will define the additive faces $F$ of the complex $\Delta\mathcal{P}$, and explain how to enumerate them knowing the additivities or non-additivities of $\verts(\Delta\mathcal{P})$ (including the limiting cones in the discontinuous case).  We will indicate the additive faces by green color and other faces by white color on the diagram. 
%\tred{make a subbadditive example and plot 2d diagram.}
See %% \autoref{fig:2d_diagrams_continuous_function} and
\autoref{fig:2d_diagrams_discontinuous_function}.
%% for illustrations of the continuous and discontinuous cases.

\textbf{For a continuous subadditive function $\pi$}, we say that a face $F \in
\Delta\mathcal{P}$ is \emph{additive} if $\Delta\pi =0$ over all $F$.  Note that
$\Delta\pi$ is affine linear over $F$, and so the face $F$ is additive if and
only if $\Delta\pi(x, y) = 0$ for all $(x, y) \in \verts(F)$. % This can be checked algorithmically.
It is clear that any subface $E$ of an additive face $F$ ($E \subseteq F$, $E
\in \Delta\mathcal{P}$) is still additive. 
Thus the additivity domain of~$\pi$ can be represented by the list of inclusion-maximal additive faces of $\Delta\mathcal{P}$; see
\cite[Lemma~3.12]{igp_survey}.\footnote{This list is computed by 
  \sage{generate\_maximal\_additive\_faces($\pi$)}.} 

\textbf{For a discontinuous subadditive function $\pi$}, we say that a face $F
\in \Delta\mathcal{P}$ is \emph{additive} if $F$ is contained in a face $F'
\in \Delta\mathcal{P}$ such that $\Delta\pi_{F'}(x,y) =0$ for any $(x,y) \in
F$.\footnote{Summarizing the detailed additivity and additivity-in-the-limit 
  situation of the function using the notion of additive faces is justified by 
  \cite[Lemmas 2.7 and 4.5]{basu-hildebrand-koeppe:equivariant} and their
  generalizations. 
  %%% NOTE: Lemma 4.5 is only stated there for the equally spaced case. 
  %%% algo-paper needs similar, more general lemmas.
}  
Since $\Delta\pi$ is affine linear in the relative interiors of each face of $\Delta\mathcal{P}$, the last condition is equivalent to $\Delta\pi_{F'}(x,y) =0$ for any $(x,y) \in \verts(F)$% , which can be checked algorithmically
. Depending on the dimension of $F$, we do the following.
\begin{enumerate}
\item Let $F$ be a two-dimensional face  of $\Delta\mathcal{P}$. If $\Delta\pi_{F}(x,y) =0$ for any $(x,y) \in \verts(F)$, then $F$ is additive. Visually on the 2d-diagram with cones, each vertex of $F$ has a green cone sitting inside $F$.
\item Let $F$ be a one-dimensional face, i.e., an edge of
  $\Delta\mathcal{P}$. Let $(x_1, y_1), (x_2, y_2)$ be its vertices. Besides
  $F$ itself, there are two other faces $F_1, F_2 \in \Delta\mathcal{P}$ that
  contain $F$. If $\Delta\pi_{F'}(x_1,y_1)=\Delta\pi_{F'}(x_2,y_2) =0$ for $F'
  = F$, $F_1$, or $F_2$, then the edge $F$ is additive. 
%%%% Visually on the 2d-diagram with cones, either of the following situations yields an additive $F$, assuming that $F$ is horizontal. 
  %%%% \tred{Illustrations are too hard to make.} \\
%%%% \hspace{5em}
%%%% \includegraphics[width=6em]{graphics-for-algo-paper/additive_edge_disct_case_1.png}
%%%% \hspace{5em}
%%%% \includegraphics[width=6em]{graphics-for-algo-paper/additive_edge_disct_case_2.png}
%%%% \hspace{5em}
%%%% \includegraphics[width=6em]{graphics-for-algo-paper/additive_edge_disct_case_3.png}
\item Let $F$ be a zero-dimensional face of $\Delta\mathcal{P}$, $F = \{(x, y)\}$. If there is a face $F' \in \Delta\mathcal{P}$ such that $(x,y) \in F'$ and $\Delta\pi_{F'}(x,y)=0$, then $F$ is additive.  Visually on the 2d-diagram with cones, the vertex $(x,y)$ is green or there is a green cone pointing at $(x,y)$. 
\end{enumerate}

On the diagram %% \autoref{fig:2d_diagrams_continuous_function} (right) and 
\autoref{fig:2d_diagrams_discontinuous_function} (right), the %(maximal)
additive faces are shaded in green.  
The projections $p_1(F)$, $p_2(F)$, and $p_3(F)$ of a two-dimensional additive face $F$ are shown as gray shadows on the $x$-, $y$- and $(x+y)$-axes of the diagram, respectively.

\section{Additional functionality}

\begin{description}
\item[\sage{minimality\_test($\pi$)}] implements a fully automatic test whether a
  given function is a minimal valid function, using the information that the
  described 2-dimensional diagrams visualize.  The algorithm is equivalent to
  the one described, in the setting of discontinuous pseudo-periodic
  superadditive functions, in Richard, Li, and Miller \cite[Theorem
  22]{Richard-Li-Miller-2009:Approximate-Liftings}.
\item[\sage{extremality\_test($\pi$)}] implements a grid-free generalization
  of the automatic extremality test from
  \cite{basu-hildebrand-koeppe:equivariant}, which is suitable also for
  piecewise linear functions with rational breakpoints that have huge
  denominators.  Its support for functions with algebraic irrational
  breakpoints such as \sagefunc{bhk_irrational}
  \cite[section~5]{basu-hildebrand-koeppe:equivariant} is experimental and
  will be described in a forthcoming paper. 
\item[\sage{generate\_covered\_intervals($\pi$)}] computes connected
  components of covered (affine imposing \cite{basu-hildebrand-koeppe:equivariant}) intervals.  This is an ingredient in
  the extremality test.
\item[\sage{extreme\_functions}] is the name of a Python module that gives access to
  the electronic compendium of extreme functions; see
  \cite{zhou:extreme-notes} and \cite[Tables 1--4]{igp_survey}.  
\item[\sage{procedures}] % is the name of a Python module that gives access to
  % the ``procedures'' that can transform 
  provides transformations of 
  extreme functions; see \cite[Table 5]{igp_survey}.
\item[\sage{random\_piecewise\_function()}] generates a random piecewise linear function with
  prescribed properties, to enable experimentation and exploration.
\item[\textsf{demo.sage}] demonstrates further functionality and the use of
  the help system.
\end{description}

%%\clearpage
\providecommand\ISBN{ISBN }
\bibliographystyle{../amsabbrvurl}
\bibliography{../bib/MLFCB_bib}

\end{document}